\documentclass[ 12pt,psamsfonts]{article}

\usepackage{theorem}
\usepackage{amscd}
\usepackage{babel}
\usepackage[leqno]{amsmath}
\usepackage{amssymb}
\usepackage{amsfonts}

\theoremstyle{change}
\theorembodyfont{\normalfont}
\newtheorem{atom}       {}      [subsection]
\newtheorem{remark}     [atom]{\textmd{\textit{Remark.}}}

\theoremheaderfont{\normalfont\bfseries\rmfamily}
\theorembodyfont{\normalfont\itshape}

\newtheorem{proposition}  [atom]{Proposition.} 

\newtheorem{theorem}      [atom]{Theorem.}

\newcommand{\ZZ}{\mathbb{Z}}

\newcommand  {\shDiv}   {\mathcal{D} \!\text{\textit{iv}}}

\newcommand  {\shI}     {\mathcal{I}}

\newcommand  {\shM}     {\mathcal{M}}

\newcommand  {\shL}     {\mathcal{L}}

\newcommand  {\shS}     {\mathcal{S}}

\newcommand  {\Ass}     {\operatorname{Ass}}

\newcommand  {\cyc}     {\operatorname{cyc}}

\newcommand  {\Div}     {\operatorname{Div}}

\renewcommand  {\k}       {\kappa}

\newcommand  {\lra}     {\longrightarrow}

\newcommand  {\primid}  {\mathfrak{p}}
\renewcommand{\O}       {\mathcal{O}}

\newcommand  {\Pic}     {\operatorname{Pic}}

\newcommand  {\ra}      {\rightarrow}

\newcommand  {\red}     {{\operatorname{red}}}

\newcommand  {\Spec}    {\operatorname{Spec}}
\newcommand  {\Supp}    {\operatorname{Supp}}

\begin{document}

\setcounter{subsection}{0}
\setlength{\unitlength}{1ex}

\begin{titlepage}
\title {Remarks on the existence of Cartier divisors}
\author{Stefan Schr\"oer\\
         }
\end   {titlepage}


\maketitle

\begin{abstract}
We characterize those invertible sheaves on a noetherian scheme
which are definable by   Cartier divisors and correct an erroneous
counterexample in the literature. 
\end{abstract}


\renewcommand{\thefootnote}{}
\footnote{\hspace{-4ex}{\itshape Key words:}  Cartier divisor,
invertible sheaf. }
\footnote{\hspace{-4ex}{\itshape Mathematics subject 
classification (1991)}:  14C20}

\subsection{Introduction}

Let 
$X $ be a noetherian scheme and 
$\shL $ an invertible 
$\O_{X} $-module; does there exist a Cartier divisor 
$D\in\Div(X) $ with 
$\shL\simeq\O_{X}(D)  $? This is no problem if 
$X $ satisfies Serre's condition 
$(S_1) $, and the   issue is to deal with embedded components.
The goal of
this short note is to provide an answer   and
to correct an erroneous counterexample in the literature.

The question was first posed by Nakai
\cite[p.\ 300]{Nakai 1963}, and later Grothendieck 
\cite[21.3.4]{Grothendieck 1967} showed that the canonical map
$\Div(X)\ra\Pic(X) $ is surjective if the subset
$\Ass(\O_{X}) \subset X$ allows an {\itshape affine}  open
neighborhood. On the other hand, it seemed to be well known from
the beginning that  in general obstructions might arise. 
Hartshorne proposed a construction
 (attributed to Kleiman)
of a non-projective irreducible 3-fold
$X $ with a {\itshape single} 
embedded component 
$x \in X $  for which it is claimed that 
$\Div(X)\ra\Pic(X) $ is   not  surjective
\cite[ex.\ 1.3, p.\ 9]{Hartshorne 1966}. Unfortunately, 
$\Ass(\O_{X})=\left\{ x ,\eta  \right\} $ is contained in every
affine open neighborhood 
$U\subset X $ of 
$x$, and Grothendieck's criterion tells us that  that
the proposed construction does {\itshape not}  yield an invertible
sheaf without Cartier divisor.

\bigskip
In the first part of this note we will discuss how the
construction can be modified in order to obtain the desired
counterexample. In the second part we will prove a
positive result, which complements Grothendieck's criterion   in
the following way: Let 
$T\subset X $ be a finite subset containing 
$\Ass(\O_{X}) $; then there is a Cartier divisor 
$D\in\Div(X) $ with 
$\shL\simeq\O_{X}(D) $ and 
support 
$\Supp(D) $ disjoint from 
$T $ if and only if the restriction of 
$\shL $ to   
$T $  is trivial.
 Here we view 
$T $ also as a ringed space, endowed with the subspace topology
and sheaf of rings 
$\O_T=i^{-1}(\O_{X}) $, where 
$i:T\ra X $ is the inclusion map.

\subsection{Absence of Cartier divisors}

In this section we construct two schemes 
$X $ for which 
$\Div(X)\ra\Pic(X) $ is not surjective.

\begin{atom}
Let us recall Hartshornes construction. We fix a ground
field 
$K $; then there is a regular, integral, proper 3-fold 
$Y $ containing two irreducible curves 
$A,B\subset Y $ such that 
$A+B $ is numerically trivial. Such a scheme is obviously
non-projective, and was constructed by Hironaka using local
blow-ups; the construction is thoroughly discussed in 
\cite[p.\ 75]{Shafarevich 1994}.  For each Cartier divisor 
$D\in \Div(Y) $ we have 
$$
A\cdot D>0 \Leftrightarrow B\cdot D<0, 
$$
and the complement of an  affine open neighborhood 
$U\subset Y $ of the generic point of
$ A $ defines such a Cartier divisor. Choose a closed point 
$a\in A $ and consider the infinitesimal
extension 
$Y\subset X $ with ideal 
$\shI=\k(a) $. The outer groups   in the exact sequence 
$$
H^1(Y,\shI) \lra \Pic(X)  \lra \Pic(Y)  \lra H^2(Y,\shI)   
$$
vanishes, hence there is an invertible 
$ \O_{X}$-module
$\shL $ with 
$B\cdot c_1(\shL)>0 $. Grothendieck's criterion tells us that 
$\shL $ is representable by a Cartier divisor
$D\in\Div(X) $; assume that 
it is even  representable by an {\itshape effective} Cartier
divisor
$D\subset X $.
But 
$A\cdot D<0 $ implies 
$A\subset D $, hence 
$a\in D $; on the other hand, according to
\cite[3.1.9]{Grothendieck 1965}, 
$D $ must be disjoint to 
$\Ass(\O_{X}) $, contradiction. 
In other words, the construction only yields  a Cartier divisor
$D\in\Div(X)$  not linearly equivalent to an effective one such
that the restriction to 
$X^\red=Y $ is equivalent to an effective Cartier divisor.

In order to achieve the desired effect we have to introduce at
least {\itshape two}  embedded components. Choose    
closed points 
$a\in A $ and 
$b\in B $, and let 
$Y\subset X $ be the infinitesimal extension with ideal 
$\shI=\k(a)\oplus \k(b) $. Again there is an invertible 
$\O_{X} $-module 
$\shL $ with 
$A\cdot c_1(\shL)<0 $ and 
$B\cdot c_1(\shL)>0 $.  We observe that  
$\Div(X)\subset Z^1(X) $ is the subgroup generated by all prime
cycles disjoint to 
$\left\{ a,b \right\} $. Assume that there is a Cartier divisor 
$D\in\Div(X) $ representing 
$\shL $. Decomposing
$D=\sum n_iD_i  $ into prime cycles, we see that each summand is
Cartier, hence
$A\cdot D_i \neq0$ and 
$B\cdot D_i \neq 0$ holds for some index 
$i $. Consequently we have   
$A\cdot D_i<0 $ and 
$a\in D_i $, or 
$B\cdot D_i<0 $ and 
$b\in D_i $; in both cases, 
$\Ass(\O_{X}) = \left\{  a,b,\eta  \right\} $ is not disjoint to 
$D_i\subset X $, contradiction. Hence it is impossible to
represent 
$\shL $ by a Cartier divisor.
\end{atom}

\begin{atom}
Another counterexample features non-separated schemes.
Let 
$A $ be a discrete valuation ring with  field of fraction 
$R $. We can glue two copies 
$U_1,U_2 $ of 
$\Spec(A) $ along 
$\Spec(R) $ and obtain an integral, regular curve 
$Y $, which is a non-separated scheme
\cite[8.8.5]{Grothendieck 1970}. The group 
$\Div(Y)=Z^1(Y) $ is isomorphic to 
$\ZZ^2 $, and the exact sequence 
$$
1\lra \Gamma(Y,\O_{Y})^\times \lra  \Gamma(Y,\shM_Y)^\times \lra 
\Div(Y)                 \lra  \Pic(Y)               \lra 0 
$$
yields 
$\Pic(Y)=\ZZ $. Let 
$Y\subset X $ be the infinitesimal extension with the ideal 
$\shI=\k(y_1)\oplus\k(y_2) $, where 
$y_1,y_2\in Y $ are the closed points. The restriction map 
$\Pic(X)\ra\Pic(Y) $ is bijective, but the sheaf 
$\shDiv_X $ is zero. Thus we have 
$\Div(X)=0 $, and 
$\shL=\O_{X}  $ is the only invertible sheaf associated to a
Cartier divisor.
\end{atom}

\subsection{Existence of Cartier divisors}

In this section, 
$X $ is a noetherian scheme, and 
$T\subset X $ is a finite subset  containing the finite subset 
$\Ass(\O_{X})\subset X $.

\begin{atom}
Let 
$\Div_T(X)\subset \Div(X) $ be the subgroup of Cartier divisors 
$D $ with 
$\Supp(D)\cap T=\emptyset $. Recall that the support 
$\Supp(D) $ is defined as the support of 
$D_1\cup D_2 $, where 
$\cyc(D)=D_1-D_2 $ is the decomposition into positive and negative
parts of the associated Weil divisor.

This construction can be sheafified: Let 
$\shS_T\subset\O_{X} $ be the subsheaf of sets whose stalk 
$\shS_{T,x} $ consists of the stalks 
$s_x\in\O_{X,x} $ whose localizations 
$s_y\in\O_{X,y} $ are units for all 
$y\in\Spec(\O_{X,x})\cap T $.  Let 
$\shM_{X,T}=\shS_T^{-1}\O_{X} $
 be the localization in the category of  sheaves of rings.
 We now define a sheaf of abelian groups 
$\shDiv_{X,T} $, written additively, by the exact sequence 
\begin{equation}
\label{sheaf of divisors}
1 \lra \O_{X}^\times  \lra \shM_{X,T}^\times  \lra  
\shDiv_{X,T} \lra 0, 
\end{equation}
and obtain 
$\Div_T(X)=\Gamma(X,\shDiv_{X,T}) $. Now let 
$i:T\ra X $ be the inclusion map, and set 
$\O_{T}=i^{-1}(\O_{X}) $. We observe the following
\end{atom}

\begin{proposition}
\label{meromorphic=local}
The 
$\O_{X} $-algebras 
$\shM_{X,T} $ and 
$i_*(\O_{T}) $ are canonically isomorphic.
\end{proposition}

First, assume that 
$X $ is the spectrum of a local ring 
$A $ with closed point 
$x\in X $. Let 
$S\subset A $ be the multiplicative subset of all 
$a\in A $ with 
$a/1\in A_\primid^\times $ for all primes 
$\primid\subset A $ corresponding to points 
$t\in T $. Clearly, 
$i_*(\O_{T})_x $ and 
$(\shM_{X,T})_x $ are canonically isomorphic to 
$S^{-1}A $. 
In the general case, consider the   diagram 
$$
\begin{CD}
i_*(\O_{T})     @>>> \prod_{x\in X} i_*(\O_{T})_x    \\
@.                   @VV\simeq V\\
\shM_{X,T}          @>>> \prod_{x\in X}  (\shM_{X,T})_x,   
\end{CD} 
$$
where the horizontal maps are the canonical inclusions.
Since the bijections 
$i_*(\O_{T})_x\ra (\shM_{X,T})_x  $ are compatible with
localization, the vertical map induces the  desired bijection  
$i_*(\O_{T}) \ra \shM_{X,T}$. QED.

\bigskip
It should be noted that these 
$\O_{X} $-algebras are in general not quasi-coherent.
From the above fact we immediately obtain the following criterion:

\begin{theorem}
\label{Cartier divisor}
An invertible 
$\O_{X} $-module 
$\shL $ is representable by a Cartier divisor 
$D\in\Div(X) $ with support disjoint from 
$T $ if and only if the restriction of 
$\shL $ to 
$T $ is   trivial in 
$\Pic(T)$.
\end{theorem}

Let 
$i:T\ra X $ be the corresponding flat morphism of ringed spaces.
The exact sequence 
(\ref{sheaf of divisors}) can be rewritten as 
$$
1 \lra \O_{X}^\times  \lra i_*i^*(\O_{X}^\times)  \lra 
\shDiv_{X,T} \lra 0, 
$$
and we obtain an exact sequence 
$$
  \Div_T(X)  \lra \Pic(X)  \lra H^1(X,i_*i^*(\O_{X}^\times)).
$$
The spectral sequence for the composition 
$\Gamma\circ i_* $ gives an inclusion 
$$
0 \lra H^1(X,i_*(\O_T^\times))  \lra   H^1(T,\O_T^\times),
$$
and we end up with the exact sequence 
$$
\Div_T(X) \lra  \Pic(X) \lra \Pic(T),
$$
which is precisely our assertion. QED.

\begin{remark}
Grothendieck's criterion can be recovered from this: Assume that 
$T\subset X $ is contained in an affine open neighborhood 
$U=\Spec(A) $. If 
$S\subset A $ is the complement of the union of all primes 
$\primid\subset A $ corresponding to points 
$x\in U\cap T $, then 
$T $ is also contained in the semi-local scheme 
$V=\Spec(S^{-1}A)  $, and 
$\Pic(X)\ra\Pic(T) $ factorizes over 
$\Pic(V) $. Since the Picard group of a semi-local ring vanishes,
each invertible 
$\O_{X} $-module 
is representable by a Cartier divisor 
$D\in\Div_T(X) $.
\end{remark}


\vspace{1em}
\noindent
        Anschrift des Autors:\\
        \\
        Mathematisches Institut\\
        Ruhr-Universit\"at\\
        44780 Bochum\\
        Germany\\ 
        E-mail s.schroeer@ruhr-uni-bochum.de

\end{document}